\documentclass[12pt,a4paper]{article}
\usepackage[a4paper, textheight=24cm,textwidth=16cm]{geometry}

\usepackage{amssymb}
\usepackage{amstext}
\usepackage{amsmath}
\usepackage{mathtools} 
\usepackage{amsthm}
\usepackage{amscd}
\usepackage{tikz-cd} \usetikzlibrary{decorations.pathmorphing}
\usepackage{xcolor} 
\usepackage{url}
\usepackage{mathrsfs}
\usepackage{hyperref}
\usepackage[utf8]{inputenc}
\usepackage[OT2,T1]{fontenc}
\newcommand{\Spec}{\mathop{\mathrm{Spec}}}
\newcommand{\Proj}{\mathop{\mathrm{Proj}}}
\newcommand{\Ker}{\mathop{\mathrm{Ker}}}

\newcommand{\Sym}{\mathop{\mathrm{Sym}}\nolimits}
\newcommand{\uSym}{\mathop{\ul{\mathrm{Sym}}}\nolimits}
\newcommand{\Supp}{\mathop{\mathrm{Supp}}}
\newcommand{\rH}{\mathrm{H}}
\newcommand{\im}{\mathop{\mathrm{Im}}}

\newcommand{\flis}{\xrightarrow{\sim}}
\newcommand{\ffl}{\longrightarrow}

\newcommand{\LL}{\Lambda}

\newcommand{\NN}{{\mathbb N}}
\newcommand{\PP}{{\mathbb P}}

\def\Fp{\mathfrak{p}}
\def\Fm{\mathfrak{m}}

\newcommand{\cK}{{\mathscr{K}}}

\newcommand{\cF}{{\mathscr{F}}}
\newcommand{\cI}{{\mathscr{I}}}
\newcommand{\cL}{{\mathscr{L}}}

\newcommand{\cO}{{\mathscr{O}}}

\newcommand{\cC}{{\mathscr{C}}}

\newcommand{\cM}{{\mathscr{M}}}

\newcounter{nc}
\renewcommand{\thenc}{{\rm(\roman{nc})}}
\newenvironment{romlist}%
{\begin{list}{\thenc}{
\usecounter{nc}
\parsep=0pt
\setlength  \labelwidth{\leftmargin}
\addtolength\labelwidth{-\labelsep}
}
}{\end{list}}
\newcounter{ctnum}
\renewcommand{\thectnum}{\textup{(\arabic{ctnum})}}
\newenvironment{numlist}%
{\begin{list}{\thectnum}{
\usecounter{ctnum} 
\parsep=0pt
\leftmargin=0pt%
\setlength{\itemindent}{\labelwidth}%
\addtolength{\itemindent}{\labelsep}%
}
}{\end{list}}

\usepackage[english]{babel}
\newtheorem{enMBteo}[subsection]{Theorem}

\newtheorem{enMBdefi}[subsection]{Definition}
\theoremstyle{remark}

\newtheorem{enMBrems}[subsection]{Remarks}

\newtheorem{enMBsubrem}[subsubsection]{Remark}

\theoremstyle{plain}

\newtheorem{enMBsubprop}[subsubsection]{Proposition}
\newtheorem{enMBsublem}[subsubsection]{Lemma}
\newtheorem{enMBsubcor}[subsubsection]{Corollary}
\newtheorem{enMBsubdefi}[subsubsection]{Definition}

\newcommand\rref[1]{{\rm\ref{#1}}}

\newcommand{\moins}{\smallsetminus}
\newcommand{\pf}{\noindent{\slshape Proof. }}
\newcommand{\ul}[1]{\underline{#1}}
\newcommand{\til}[1]{\widetilde{#1}}
\newcommand{\Tag}[1]%
{\cite[\href{https://stacks.math.columbia.edu/tag/#1}{Tag #1}]{StProj}}
\title{The relative Fujita-Zariski theorem}
\author{%
Laurent Moret-Bailly\thanks{Univ Rennes, CNRS, IRMAR - UMR 6625, F-35000 Rennes, France}
 \thanks{\tt laurent.moret-bailly[AT]univ-rennes.fr}
}

\begin{document}
\selectlanguage{english}
\date{\today}
\maketitle
%
\selectlanguage{english}
\begin{abstract}
We prove, with no claim to originality, a relative version of the Fujita-Zariski theorem. When the base is a field, this result is  due to Fujita \cite{Fu} and states that if an invertible sheaf on  a proper variety is ample on its base locus, its sufficiently high powers are globally generated. The special case where the base locus is finite was proved by Zariski \cite{Zar}, whence the name.
\end{abstract}

\setcounter{tocdepth}{2}
\tableofcontents

\section{Definitions, notation, and statement}
Throughout the paper we consider a noetherian ring $R$ and a proper morphism of schemes $f:X\to S:=\Spec(R)$. 
\begin{enMBdefi}\label{DefSysLin}
A \emph{linear system} $\Lambda=(\cL,V,\rho)$ on $X$ consists of:
\begin{romlist}
\item an invertible $\cO_{X}$-module $\cL$,
\item an $R$-module $V$ of finite type, and 
\item an $R$-linear map $\rho: V\to \rH^0(X,\cL)$, or (equivalently) an $\cO_{X}$-linear map $\til{\rho}: f^*V\to \cL$.
\end{romlist}
\end{enMBdefi}
From such a $\LL$ we derive an exact sequence $f^*V\otimes\cL^{-1}\to\cO_{X}\to\cO_{B}\to0$ where $B$ is a closed subscheme of $X$ called the \emph{base locus} of the system.

\begin{enMBteo}\label{ThFZ}
Let $\LL$ be a linear system on $X$ as above. Assume that the restriction of $\cL$ to $B$ is ample. Then there exists an integer  $t_{0}$ such that $\cL^{\otimes t}$ is globally generated for all $t\geq t_{0}$.
\end{enMBteo}
\begin{enMBrems} 
\begin{numlist}
\item Note that the conclusion of the theorem is about $\cL$ only. The assumption is preserved if we replace $V$ by its image in $\rH^0(X,\cL)$ or a bigger submodule (since this can only shrink $B$). So we could as well assume that $V=\rH^0(X,\cL)$. However, this does not make the proof any  simpler. Moreover, given an element $v\in V$, it is advisable to distinguish $v$ (living in an $R$-module) from $\rho(v)$ which is a section of $\cL$ and thus lives on $X$.
\item When $R$ is a field, this result is due to Fujita \cite[Theorem 1.10]{Fu}. The special case when $B$ is finite was proved earlier by Zariski \cite[Theorem 6.2]{Zar}. 
\item Our proof follows Fujita's strategy rather faithfully. The only difficulties were to find the right substitutes for induction on dimension of a subscheme, on dimension of $V$, and the notion of a ``general element'' of $V$. 
\item The interested reader will have no trouble adapting the proof to the case where $X$ is a proper algebraic space over $S$.
\item To the author's knowledge, and to this date, Theorem \ref{ThFZ} is not in the literature. It is used without proof in the paper \cite{Sch}, with a reference to Fujita's article.
\end{numlist}
\end{enMBrems}
\subsection{Graded rings and modules}
By convention, all our graded rings and modules are $\NN$-graded. For modules, this may seem somewhat unnatural, but it is harmless as the properties we shall be concerned with can be checked ``in large enough degree''.

Let  $X\to S$ and $\LL$ be as above. To a coherent $\cO_{X}$-module $\cF$ and an integer $q\geq0$ we associate the graded $R$-module
$$\rH^q_{*}(\cF)=\rH^q(X,\cF\otimes\uSym_{\cO_{X}}(\cL))=\bigoplus_{t\geq0}\rH^q(X,\cF\otimes\cL^{\otimes t})$$
whose graded components are finite $R$-modules. Note that $\rH^0_{*}(\cO_{X})=\rH^0(X,\uSym_{\cO_{X}}(\cL))$ is a graded $R$-algebra, and each $\rH^q_{*}(\cF)$ is a graded $\rH^0_{*}(\cO_{X})$-module. Moreover, $\rho:V\to \rH^0(X,\cL)$ gives rise to a morphism 
$$\Sym_{R}(V)\to \rH^0_{*}(\cO_{X})$$
of graded $R$-algebras, thus turning $\rH^q_{*}(\cF)$ into a graded $\Sym_{R}(V)$-module. Technically, these modules will be our main object of study; observe that $\Sym_{R}(V)$ is a finitely generated $R$-algebra and thus a noetherian ring, unlike $\rH^0_{*}(\cO_{X})$ in general.

We shall make use of the following easy fact: a graded  $\Sym_{R}(V)$-module $M=\bigoplus_{t\geq0}M_{t}$ is finitely generated if and only if each $M_{t}$ is finitely generated over $R$, and the natural map $V\otimes_{R}M_{t}\to M_{t+1}$ (where $V$ is viewed as the degree 1 component of $\Sym_{R}(V)$) is surjective for $t$ large enough.
\begin{enMBsubdefi}\label{DefModSys} A \emph{$\LL$-module} $(\cF,q,M_{*})$ consists of
\begin{romlist}
\item a coherent $\cO_{X}$-module $\cF$,
\item an integer $q\geq0$,
\item a graded $\Sym_{R}(V)$-submodule $M_{*}=\bigoplus_{t\geq0}M_{t}\subset \rH^q_{*}(\cF)$.
\end{romlist}

\end{enMBsubdefi}
This differs slightly from the corresponding notion of a ``$\LL$-module system'' in \cite{Fu}, where negative degrees are allowed.
Such a triple $(\cF,q,M_{*})$ will usually be called $M_{*}$ for simplicity. In particular, we say that $(\cF,q,M_{*})$  is \emph{finitely generated} if $M_{*}$ is a finitely generated $\Sym_{R}(V)$-module. 

\begin{enMBsubrem} 
To prove Theorem \ref{ThFZ}, an essential tool is to prove that (under its assumptions) certain  $\LL$-modules are finitely generated. We shall use  induction arguments typically involving the choice of a ``general'' element of $V$. The next section explains what we mean by this.

\end{enMBsubrem}
\subsection{Fat subsets in modules, general position results}\label{Fat}
\begin{enMBsubdefi}\label{DefFat} Let $A$ be a ring, $M$ an $A$-module. A subset $E\subset M$ is \emph{($A$-)fat} if $M\moins E$ is contained in a finite union of proper submodules of $M$.
\end{enMBsubdefi}
\begin{enMBsublem}\label{LemFat} Let $A$ be a ring. 
\begin{numlist}
\item\label{LemFat1}The following are equivalent:
\begin{romlist}
\item\label{LemFat11} For every $A$-module $M$, every $A$-fat subset of $M$ is nonempty.
\item\label{LemFat12} For every finitely generated $A$-module $M$, every $A$-fat subset of $M$ is nonempty.
\item\label{LemFat13} For each maximal ideal $\Fm$ of $A$, the residue field $A/\Fm$ is infinite.
\item\label{LemFat14} For each prime ideal $\Fp$ of $A$, the residue field $\kappa(\Fp)$ is infinite.
\end{romlist}
\item\label{LemFat2} There is a faithfully flat $A$-algebra $A'$ satisfying the conditions of \rref{LemFat1}. Moreover we can take $A'$ noetherian if $A$ is.
\end{numlist}
\end{enMBsublem}

\pf  \ref{LemFat1} The implication \ref{LemFat11}$\Rightarrow$\ref{LemFat12} is trivial. For the converse, let $M$ be an $A$-module and $(N_{i})_{1\leq i\leq r}$ a finite family of submodules. Pick an element in $M\moins N_{i}$ for each $i$, and apply \ref{LemFat12} to the submodule generated by these. 

\ref{LemFat12}$\Rightarrow$\ref{LemFat13}: Assume $A$ has a finite residue field $\kappa$. Then the finitely generated $A$-module $\kappa^2$ is the union of its finitely many proper submodules.

\ref{LemFat13}$\Leftrightarrow$\ref{LemFat14} is immediate. Now we assume \ref{LemFat13} and prove \ref{LemFat12}. Let $M$ and $(N_{i})_{1\leq i\leq r}$ be as above, with $M$ finitely generated, and let us prove that  $\bigcup_{1\leq i\leq r}N_{i}\varsubsetneq M$.
For each $i$, $M/N_{i}$ is  finitely generated, so there is a maximal ideal $\Fm_{i}$ and an epimorphism $\pi_{i}:M/N_{i}\twoheadrightarrow \kappa_{i}:=A/\Fm_{i}$. Clearly we may replace $N_{i}$ by $\Ker\pi_{i}$. In other words, $N_{i}$ is now the preimage of a hyperplane $H_{i}$ in the $\kappa_{i}$-vector space $M\otimes_{A}\kappa_{i}$. 

First assume that all the $\Fm_{i}$'s are equal, with residue field $\kappa$. Since $\kappa$ is infinite by assumption, we have $\bigcup_{1\leq i\leq r}H_{i}\varsubsetneq M\otimes_{A}\kappa$, whence the result.

Otherwise we renumber the $\kappa_{i}$'s as $\kappa_{1},\dots,\kappa_{s}$ (assumed pairwise distinct), and the $N_{i}$'s as $N_{jk}$ in such a way that $M/N_{jk}\cong\kappa_{j}$. For each $j$ we pick an element $m_{j}$ of $M\otimes\kappa_{j}$ which is not in the image of any $N_{jk}$, as in the previous step. Then we observe that the $\Fm_{j}$'s are pairwise coprime, and therefore $M\to \prod_{j}M\otimes\kappa_{j}$ is surjective, so there is $m\in M$ which reduces to $m_{j}$ for all $j$ and thus cannot belong to any $N_{jk}$.\medskip

\noindent\ref{LemFat2} We can take $A'=U^{-1}A[T]$ where $U$ is the multiplicative set of all monic polynomials. (If $A$ is noetherian, another choice is $A'=A[[T]][T^{-1}]$).
\qed
\begin{enMBsubprop}\label{PropPosGen} Let $R$, $f:X\to S$ and $\LL=(\cL,V,\rho)$ be as in \rref{DefSysLin}, with base locus $B\subset X$. Let $\cF$ be a coherent $\cO_{X}$-module, and let $\Sigma\subset X\moins B$ be a finite set of points. Then, for all $\delta$ in a suitable fat subset of $V$, the following conditions are satisfied:
\begin{romlist}
\item\label{PropPosGen1} The section $\rho(\delta)$ of $\cL$ does not vanish at any  $x\in\Sigma$.
\item\label{PropPosGen2} The morphism $\varphi:\cF\otimes\cL^{-1}\to\cF$ induced by tensoring with $\rho(\delta)$ is injective on $X\moins B$; in other words, $\Supp\,(\Ker\varphi)\subset B$.
\end{romlist}
\end{enMBsubprop}
\pf We can enlarge $\Sigma$ and assume it contains $\mathop{\mathrm{Ass}_{\cO_{X}}}(\cF)\moins B$. Then condition \ref{PropPosGen2} is a consequence of  \ref{PropPosGen1}. Now for $x\in\Sigma$, $V_{x}:=\left\{\delta\in V\mid \rho(\delta)(x)=0\right\}$ is an $R$-submodule of $V$, and $V_{x}\varsubsetneq V$ since $x\notin B$. The result follows.\qed
\medskip

Here is a first consequence:
\begin{enMBsubcor}\label{LemModGrad} Let $R$ be a noetherian  ring, $V$ a finitely generated $R$-module, $M_{*}$  a finitely generated graded $\Sym(V)$-module. There exists $t_{0}\geq0$ and a fat subset $E$ of $V$ such that for all $\delta\in E$, the map $M_{t}\xrightarrow{\times\delta}M_{t+1}$ is injective for all $t\geq t_{0}$.
\end{enMBsubcor}
\pf We apply Proposition \ref{PropPosGen} to the following data: let $P:=\PP(V)=\Proj(\Sym(V))$, $\cM=\til{M_{*}}$ the coherent $\cO_{P}$-module associated to $M_{*}$. We take for $\cL$ the canonical sheaf $\cO_{P}(1)$; note that the base locus $B$ is empty. By the proposition, for $\delta$ in a fat subset of $V$ the corresponding $\cM(-1)\to \cM$ is injective, so of course the same holds for all maps $\cM(t)\to \cM(t+1)$ and $\rH^0(P,\cM(t))\to \rH^0(P,\cM(t+1))$. 

On the other hand, by \cite[(2.3.1)]{EGA2} there exists $t_{0}$ such that for all $t\geq t_{0}$ the  canonical map $M_{t}\to \rH^0(P,\cM(t))$ is an isomorphism, which completes the proof.\qed

\section{Proof of the theorem}
\subsection{A finiteness result}
A key step in the proof is the following proposition (which will be used for $q=1$ only!):
\begin{enMBsubprop}\label{PropTF} With the assumptions of Theorem \rref{ThFZ}, let $(\cF,q,M_{*})$ be a $\LL$-module with $q>0$. Then $M_{*}$ is finitely generated.
\end{enMBsubprop}

The proof relies on a d\'evissage lemma:
\begin{enMBsublem}\label{LemDeviss} Let $\LL=(\cL,V,\rho)$ be a linear system on $X$. Let $\delta\in V$ be  fixed, and let $\varphi:\cF\otimes\cL^{-1}\to\cF$ denote  the morphism given by multiplication by $\rho(\delta)$. Consider the self-defining exact sequence
$$0\ffl \cK\ffl \cF\otimes\cL^{-1}\xrightarrow{\;\;\varphi\;\;}\cF\xrightarrow{\;\;\pi\;\;} \cC\ffl0$$
of sheaves on $X$. We fix an integer $q\geq0$, and we introduce the following $\LL$-modules:
$$\begin{array}{rcl}
M_{*} & := & \rH^q_{*}(\cF)\\
N_{*} & := & \im \left(M_{*}\to \rH^q_{*}(\cC)\right)\\
K^+_{*} & := & \rH^{q+1}_{*}(\cK)
\end{array}$$
and we assume the following conditions:
\begin{romlist}
\item\label{LemDeviss1} $K^+_{t}=0$ for $t\gg0$;
\item\label{LemDeviss2} $N_{*}$ is a finitely generated $\LL$-module.
\end{romlist}
Then $M_{*}$ is finitely generated.
\end{enMBsublem}
\noindent{\slshape Proof of \ref{LemDeviss}.} Condition \ref{LemDeviss2} means that $VN_{t-1}=N_{t}$ for large $t$. We fix $t$ such that this holds in addition to $K^+_{t}=0$, and we proceed to show that $VM_{t-1}=M_{t}$, which will prove the result.

To achieve this we already observe that we have a natural surjective morphism $M_{*}\twoheadrightarrow N_{*}$ of graded $\Sym(V)$-modules. This and the condition  $VN_{t-1}=N_{t}$  imply that $VM_{t-1}\hookrightarrow M_{t}\to N_{t}$ is surjective. Thus it suffices to prove that $\Ker(M_{t}\to N_{t})\subset VM_{t-1}$. We shall see that, more precisely, $\Ker(M_{t}\to N_{t})=\delta M_{t-1}$.

Le $\cI\subset\cF$ be the image of $\varphi$. Put $I_{*}=\rH^q_{*}(\cI)$. The map $M_{t-1}\to M_{t}$ induced by $\varphi$ is just multiplication by $\delta$ for the $\Sym(V)$-module structure on $M_{*}$ (in particular its image is $\delta M_{t-1}$), and it factors as $M_{t-1}\xrightarrow{\alpha}I_{t}\xrightarrow{\beta}M_{t}$. From the short exact sequence $0\to\cK\to\cF\otimes\cL^{-1}\to\cF\to \cI\ffl0$ (twisted by $\cL^{\otimes t}$) and the condition on $K^+_{t}$ we infer that $\alpha:M_{t-1}\to I_{t}$ is surjective. Thus $\im\beta=\im(\beta\circ\alpha)=\delta \,M_{t-1}$. On the other hand,  the short exact sequence $0\to\cI\to\cF\to\cC\to0$ yields an exact sequence $I_{t}\xrightarrow{\beta} M_{t}\to N_{t}\to0$. Combining these we conclude that 
$\Ker(M_{t}\to N_{t})=\delta\,M_{t-1}$, as promised.\qed\medskip

\noindent{\slshape Proof of \ref{PropTF}.}  We first note that the question is local on $S$ for the fpqc topology. In particular, applying Lemma \ref{LemFat}\,\ref{LemFat2}, we may assume that all the residue fields of $R$ are infinite, so that fat subsets of $R$-modules are always nonempty.

Next, we may and will assume that $M_{*}$ is the full $\rH^q_{*}(\cF)$ since $\Sym(V)$ is a noetherian ring.

We now choose a point $s\in S$ and work locally, i.e.\ we shall find an affine neighborhood $\Spec(R')$ of $s$ such that $M_{*}\otimes_{R}R'$ is a finitely generated $\Sym(V)\otimes_{R}R'$-module. 
Putting $Y:=\Supp(\cF)$ (defined, as a subscheme of $X$, by the ideal $\mathop{\mathrm{Ann}_{\cO_{X}}}(\cF)$), we proceed by induction on 
$$d(\cF):=\dim_{\kappa(s)}\im\left(V\to \rH^0(Y_{s},\cL_{\mid Y_{s}})\right).$$

If $d(\cF)=0$, then $Y_{s}\subset B$. Therefore $\cL$ is ample on $Y_{s}$. By \cite[(4.7.1)]{EGA3_I} we may assume that $\cL$ is  ample on $Y$ by restricting to a neighborhood of $s$. But then  the result is trivial since $M_{t}=0$ for large $t$. (The condition $q>0$ is used here).

Now assume $d(\cF)>0$, and the result proved for all sheaves with smaller $d$. We have $Y_{s}\not\subset B$. Fix a point $y\in Y_{s}\moins B$ and apply Proposition \ref{PropPosGen} with $\Sigma=\{y\}$. With our assumption on the residue fields, we see that there exists $\delta\in V$ such that
\begin{romlist}
\item\label{Cond1} $\rho(\delta)\in\rH^0(X,\cL)$ does not vanish identically on $Y_{s}$, and
\item\label{Cond2} $\varphi:\cF\otimes\cL^{-1}\to \cF$ given by $\rho(\delta)$ is injective on $X\moins B$.
\end{romlist}
We now apply Lemma \ref{LemDeviss}. We have an exact sequence 
$$0\ffl \cK\ffl \cF\otimes\cL^{-1}\xrightarrow{\;\;\varphi\;\;}\cF\xrightarrow{\;\;\pi\;\;} \cC\ffl0$$
where $\Supp\cK\subset B$ because of \ref{Cond2}. In particular, $\cL$ is ample on $\Supp\cK$, whence $\rH^{q+1}(\cK\otimes\cL^{\otimes t})=0$ for large $t$, which is condition  \ref{LemDeviss1} of  \ref{LemDeviss}.

Consider $N_{*}:=\im(M_{*}\to \rH^q_{*}(\cC))$. We have a $\LL$-module $(\cC,q,N_{*})$, and we shall use the induction hypothesis to prove that $N_{*}$ is finitely generated: this will complete the proof by \ref{LemDeviss}. For this, it suffices to prove that $d(\cC)<d(\cF)$. Putting $Z=\Supp(\cC)$, we have a surjection of $\kappa(s)$-vector spaces
$$\im\left(V\to  \rH^0(Y_{s},L_{\mid Y_{s}})\right)\twoheadrightarrow \im\left(V\to   \rH^0(Z_{s},L_{\mid Z_{s}})\right)$$
of dimensions $d(\cF)$ and $d(\cC)$ respectively. The image of $\delta\in V$ in the first space is nonzero by the above condition \ref{Cond1}, but it vanishes in the second space by definition of $\cC$. Thus, $d(\cC)<d(\cF)$, as claimed.\qed
\medskip

\subsection{Proof of Theorem \ref{ThFZ}}
We keep the notation and assumptions of \ref{ThFZ}. As in the proof of \ref{LemDeviss}, we assume that all residue fields of $R$ are infinite. In addition, we may assume $V=\rH^0(X,\cL)$: this does not change the conclusion, and can only make the base locus smaller.

We procced by noetherian induction on $X$, assuming that for all proper closed subschemes $Y\varsubsetneq X$, $\cL_{\mid Y}^{\otimes t}$ is globally generated for large $t$. We also assume $B\varsubsetneq X$ set-theoretically: otherwise, $\cL$ is ample. In particular, $V\neq0$.

We shall apply the induction to the subscheme $Y\varsubsetneq X$ defined by $\rho(\delta)$, for a suitable $0\neq\delta\in V$. The idea is to find $\delta$ such $\rH^0(\cL^{\otimes t})\ffl \rH^0(\cL_{\mid Y}^{\otimes t})$ is surjective for large $t$. This clearly suffices because then $\cL^{\otimes t}$ will have  no base points on $Y$ (by induction) and no base points off $Y$ since $Y$ clearly contains $B$. The surjectivity requirement translates into an injectivity property on $\rH^1$, which motivates the choices below.\medskip

\noindent I claim that there is a fat $E\subset V$ (not containing $0$) such that, for all $\delta\in E$:
\begin{romlist}
\item\label{Condd1} multiplication by $\delta$ on  $\rH^1_{*}(\cO_{X})$ is injective on components of large degree;
\item\label{Condd2} $\rho(\delta)$ is regular off $B$, i.e.\ $\varphi:\cL^{-1}\to \cO_{X}$ given by $\rho(\delta)$ is injective on $X\moins B$ .
\end{romlist}
Indeed, \ref{Condd2} is the same as in the proof of \ref{PropTF}, and for \ref{Condd1} we see by \ref{PropTF} that $\rH^1_{*}(\cO_{X})$ is a finitely generated $\Sym(V)$-module, so Corollary \ref{LemModGrad} applies. Putting $\cI=\im(\varphi)\subset\cO_{X}$, we have exact sequences (where $Y\varsubsetneq X$ is the zero locus of $\delta$)
$$\begin{array}{rccccccccl}
 0 & \ffl & \cK & \ffl & \cL^{-1} &  \ffl  & \cI &\ffl &  0\\
 0 & \ffl & \cI & \ffl & \cO_{X} &  \ffl & \cO_{Y} &\ffl &  0. 
 \end{array}$$
and the same sequences twisted by $\cL^{\otimes t}$ for all $t$. 
By \ref{Condd1}, $\cK$ has support in $B$ on which $\cL$ is ample, so the first sequence induces $\rH^1(\cL^{\otimes(t-1)})\flis\rH^1(\cI\otimes\cL^{\otimes t})$ for large $t$. The composition $\rH^1(\cL^{\otimes(t-1)})\flis\rH^1(\cI\otimes\cL^{\otimes t})\ffl\rH^1(\cL^{\otimes t})$ is given by  multiplication by $\delta$  on $\rH^1_{*}(\cO_{X})$, thus injective for large $t$ by \ref{Condd1}. In short, for $t\gg0$ the second exact sequence gives rise to injections $\rH^1(\cI\otimes\cL^{\otimes t})\ffl\rH^1(\cL^{\otimes t})$, hence the restriction maps $\rH^0(\cL^{\otimes t})\ffl \rH^0(\cL_{\mid Y}^{\otimes t})$ are surjective. As explained above, this completes the proof.\qed

\bibliographystyle{plain}
 \bibliography{Fujizar}

\begin{thebibliography}{1}

\bibitem{Fu}
Takao Fujita.
\newblock Semipositive line bundles.
\newblock {\em J. Fac. Sci. Univ. Tokyo Sect. IA Math.}, 30(2):353--378, 1983.

\bibitem{EGA2}
A.~Grothendieck.
\newblock \'{E}l\'ements de g\'eom\'etrie alg\'ebrique. {II}. \'{E}tude globale
  \'el\'ementaire de quelques classes de morphismes.
\newblock {\em Inst. Hautes \'Etudes Sci. Publ. Math.}, (8):5--222, 1961.

\bibitem{EGA3_I}
A.~{Grothendieck}.
\newblock {\'El\'ements de g\'eom\'etrie alg\'ebrique. III: \'Etude
  cohomologique des faisceaux coh\'erents. (Premi\`ere partie)}.
\newblock {\em {Publ. Math., Inst. Hautes \'Etud. Sci.}}, 11:81--167, 1961.

\bibitem{Sch}
Stefan Schr\"oer.
\newblock A characterization of semiampleness and contractions of relative
  curves.
\newblock {\em Kodai Math. J.}, 24(2):207--213, 2001.

\bibitem{Zar}
Oscar Zariski.
\newblock The theorem of {R}iemann-{R}och for high multiples of an effective
  divisor on an algebraic surface.
\newblock {\em Ann. of Math. (2)}, 76:560--615, 1962.

\end{thebibliography}

\end{document}